\documentclass[10pt]{amsart}
\usepackage{amsmath, amsfonts, amssymb, amstext, cancel, amscd, graphicx}
\usepackage[all,2cell,ps]{xy}

\begin{document}

\title[Multivariable Calculus and Differential Forms]{A Tasty Combination: \\ Multivariable Calculus and Differential Forms}

\author{Edray Herber Goins}
\address{Department of Mathematics \\ Purdue University \\ 150 North University Street \\ West Lafayette, IN 47907}
\email{egoins@math.purdue.edu}

\author{Talitha M. Washington}
\address{Department of Mathematics \\ 1800 Lincoln Avenue \\ University of Evansville \\ Evansville, IN 47722}
\email{tw65@evansville.edu}

\begin{abstract}
Differential Calculus is a staple of the college mathematics major's diet.  Eventually one becomes tired of the same routine, and wishes for a more diverse meal.  The college math major may seek to generalize applications of the derivative that involve functions of more than one variable, and thus enjoy a course on Multivariate Calculus.  We serve this article as a culinary guide to differentiating and integrating functions of more than one variable -- using differential forms which are the basis for de Rham Cohomology.
\end{abstract}

\maketitle

\section*{Introduction}

Differential Calculus is a staple of the college mathematics major's diet.  It is relatively easy to explain the Fundamental Theorem of Calculus: Given a differentiable function $f: I \to \mathbb R$ defined on closed interval $I = [a,b]$, there is the identity
\[ \int_a^b \frac d{dt} \bigl[ f(t) \bigr] dt = f(b) - f(a). \]
\noindent Eventually one becomes tired of the same routine, and wishes for a more diverse meal.  The college math major may seek to generalize applications of the derivative that involve functions of more than one variable, and thus enjoy a course on Multivariate Calculus.  Actually, there is a veritable buffet of ways to differentiate and integrate a function of more than one variable: there is the gradient, curl, divergence, path integrals, surface integrals, and volume integrals.  Plus, there are many ``Fundamental'' Theorems of Multivariate Calculus, such as Stokes' Theorem, Green's Theorem, and Gauss' Theorem.

We serve this article as a culinary guide to differentiating and integrating functions of more than one variable -- using differential forms which are the basis for de Rham Cohomology.

\section*{Gradient, Curl, and Divergence}

Let's focus on functions of three variables.  First, let's fix a ``simply connected'' closed subset $D \subseteq \mathbb R^3$; this means we can define an integral between two points without concern of a choice of path in $D$ which connects them.  (For example, a subset in the form $D = [a,b] \times [c,d] \times [p,q]$ is simply connected.)  Recall that the \emph{gradient} of a scalar function $f: D \to \mathbb R$ is the vector-valued function $\nabla f: D \to \mathbb R^3$ defined by
\[ \nabla f = \frac {\partial f}{\partial x} \, \mathbf{i} + \frac {\partial f}{\partial y} \, \mathbf{j} + \frac {\partial f}{\partial z} \, \mathbf{k}. \]
\noindent (Here $\mathbf{i}=(1,0,0)$, $\mathbf{j}=(0,1,0)$, and $\mathbf{k}=(0,0,1)$ are the standard basis vectors for 3-dimensional space.) The \emph{curl} of a vector field $\mathbf{F}: D \to \mathbb R^3$, which we write in the form $\mathbf{F} = M \, \mathbf{i} + N \, \mathbf{j} + P \, \mathbf{k}$, is the vector-valued function $\nabla \times \mathbf{F}: D \to \mathbb R^3$ defined by
\[ \nabla \times \mathbf{F} = \left| \begin{matrix} \mathbf{i} & \mathbf{j} & \mathbf{k} \\[8pt] \dfrac {\partial}{\partial x} & \dfrac {\partial}{\partial y} & \dfrac {\partial}{\partial z} \\[8pt] M & N & P \end{matrix} \right| = \left( \frac {\partial P}{\partial y} - \frac {\partial N}{\partial z} \right) \mathbf{i} + \left( \frac {\partial M}{\partial z} - \frac {\partial P}{\partial x} \right) \mathbf{j} + \left( \frac {\partial N}{\partial x} - \frac {\partial M}{\partial y} \right) \mathbf{k}. \]
\noindent (If you don't remember how to compute $3 \times 3$ determinants, don't worry; we have a short exposition on them contained in the appendix.)  Finally, the \emph{divergence} of a vector field $\mathbf{G}: D \to \mathbb R^3$, which we write in the form $\mathbf{G} = S \, \mathbf{i} + T \, \mathbf{j} + U \, \mathbf{k}$, is the scalar function $\nabla \cdot \mathbf{G}: D \to \mathbb R$ defined by
\[ \nabla \cdot \mathbf{G} = \frac {\partial S}{\partial x} + \frac {\partial T}{\partial y} + \frac {\partial U}{\partial z}. \]
\noindent Larson \cite{Larson2007Calculus:-Early} provides more information on Multivariable Calculus.  Marsden and Tromba \cite{Marsden2004Vector-Calculus} is also an excellent source.  Hughes-Hallet et al. \cite{Hughes-Hallett2004Calculus:-Singl} provides a novel approach to the subject using concept-based learning. \\

We would like to answer the following questions:
\textit{\begin{enumerate}
\item Which $\mathbf{F}: D \to \mathbb R^3$ have a scalar potential $f: D \to \mathbb R$ such that $\mathbf{F} = \nabla f$?
\item Which $\mathbf{G}: D \to \mathbb R^3$ have a vector potential $\mathbf{F}: D \to \mathbb R^3$ such that $\mathbf{G} = \nabla \times \mathbf{F}$?
\item Which $f: D \to \mathbb R$ have a vector potential $\mathbf{G}: D \to \mathbb R^3$ such that $f = \nabla \cdot \mathbf{G}$?
\end{enumerate}}

\subsection*{Relating Gradients and Curls}

Let's answer the first of our motivating questions.  We'll show the following:

\
\begin{center} \textit{$\mathbf{F} = \nabla f$ is a gradient if and only if the curl $\nabla \times \mathbf{F} = \mathbf{0}$.} \end{center}
\

\noindent First assume that $\mathbf{F} = \nabla f$.  That is, $\mathbf{F}= M \, \mathbf{i} + N \, \mathbf{j} + P \, \mathbf{k} = \dfrac {\partial f}{\partial x} \, \mathbf{i} + \dfrac {\partial f}{\partial y} \, \mathbf{j} + \dfrac {\partial f}{\partial z} \, \mathbf{k}.$  To find the curl of $\mathbf{F}$, we have
\begin{align*}
\nabla \times \mathbf{F} & =  \left( \frac {\partial P}{\partial y} - \frac {\partial N}{\partial z} \right) \mathbf{i} + \left( \frac {\partial M}{\partial z} - \frac {\partial P}{\partial x} \right) \mathbf{j} + \left( \frac {\partial N}{\partial x} - \frac {\partial M}{\partial y} \right) \mathbf{k}\\
& = \left(\frac {\partial^2 f}{\partial y \, \partial z}  - \frac {\partial^2 f}{\partial z \, \partial y}\right)\mathbf{i}+\left(\frac {\partial^2 f}{\partial z \, \partial x} - \frac {\partial^2 f}{\partial x \, \partial z}\right)\mathbf{j}+\left(\frac {\partial^2 f}{\partial x \, \partial y} - \frac {\partial^2 f}{\partial y \, \partial x}\right)\mathbf{k}\\
& = \mathbf{0}.
\end{align*}

\noindent Note that we can change the order of differentiation as long as $f$ is continuously twice-differentiable: with this assumption, Clairaut's Theorem states that the mixed partial derivatives are equal.  (This theorem and its proof can be found in \cite[pgs. 885 and A-48]{Stewart2008Calculus:-Early}.)

Now what if $\nabla \times \mathbf{F} = \mathbf{0}$?  Can we find an $f$ such that $\mathbf{F} = \nabla f$?  Well, if $\nabla \times \mathbf{F} = \mathbf{0}$, then
\[  \left( \frac {\partial P}{\partial y} - \frac {\partial N}{\partial z} \right) = \left( \frac {\partial M}{\partial z} - \frac {\partial P}{\partial x} \right) = \left( \frac {\partial N}{\partial x} - \frac {\partial M}{\partial y} \right) = 0. \]
\noindent  Now define the scalar function $f = u + v + w$ in terms of the definite integrals
\[ \begin{aligned} u(x,y,z) & = \int_{x_0}^x M(\sigma, \, y, \, z) \, d \sigma, \\[5pt] v(x,y,z) & = \int_{y_0}^y \left[ N(x, \tau, z) - \frac {\partial u}{\partial y}(x, \tau, z) \right] d \tau, \\[5pt] w(x,y,z) & = \int_{z_0}^z \left[ P(x, y, \zeta) - \frac {\partial u}{\partial z}(x, y, \zeta) - \frac {\partial v}{\partial z}(x, y, \zeta) \right] d \zeta \end{aligned} \] \noindent  where we have fixed a point $x_0 \, \mathbf{i} + y_0 \, \mathbf{j} + z_0 \, \mathbf{k} \in D$.  (Here's where we subtly use the assumption that $D$ is simply connected:  the integrals are independent of a choice of path in $D$ which connects $x_0 \, \mathbf{i} + y_0 \, \mathbf{j} + z_0 \, \mathbf{k}$ and $x \, \mathbf{i} + y \, \mathbf{j} + z \, \mathbf{k}$.)  Upon changing the order of differentiation and integration, we have the derivatives \[ \begin{aligned} \frac {\partial v}{\partial x} & = \int_{y_0}^y \left[ \frac {\partial N}{\partial x} - \frac {\partial M}{\partial y} \right] d \tau = 0, \\[5pt] \frac {\partial w}{\partial x} & = \int_{z_0}^z \left[ \left( \frac {\partial P}{\partial x} - \frac {\partial M}{\partial z} \right) - \frac {\partial v}{\partial x} \right] d \zeta = 0, \\[5pt] \frac {\partial w}{\partial y} & = \int_{z_0}^z \left[ \frac {\partial P}{\partial y} - \frac {\partial}{\partial z} \left( \frac {\partial u}{\partial y} + \frac {\partial v}{\partial y} \right) \right] d \zeta = \int_{z_0}^z \left[ \frac {\partial P}{\partial y} - \frac {\partial N}{\partial z} \right] d \zeta = 0. \end{aligned} \]
\noindent (We can only do this if $M$, $N$, and $P$ are continuously differentiable functions; this also follows from Clairaut's Theorem.) It follows that $\mathbf{F} = \nabla f$:
\[ \nabla f = \biggl[ \frac {\partial u}{\partial x} + \cancelto{0}{ \frac {\partial v}{\partial x} } + \cancelto{0}{ \frac {\partial w}{\partial x} } \biggr] \mathbf{i} + \biggl[ \frac {\partial u}{\partial y} + \frac {\partial v}{\partial y} + \cancelto{0}{ \frac {\partial w}{\partial y} } \biggr] \mathbf{j} + \biggl[ \frac {\partial u}{\partial z} + \frac {\partial v}{\partial z} + \frac {\partial w}{\partial z} \biggr] \mathbf{k}  = M \, \mathbf{i} + N \, \mathbf{j} + P \, \mathbf{k}. \]

We conclude that $\mathbf{F}: D \to \mathbb R^3$ has a scalar potential $f: D \to \mathbb R$ such that $\mathbf{F} = \nabla f$ if and only if the curl $\nabla \times \mathbf{F} = \mathbf{0}$.

\subsection*{Example}

\begin{figure}[t] \begin{center} \caption{Plot of $\mathbf{F} = \nabla f$ for $f(x,y,z) = \sqrt{x^2 + y^2 + z^2}$} \label{vector_field_1} \includegraphics[width=0.85\textwidth]{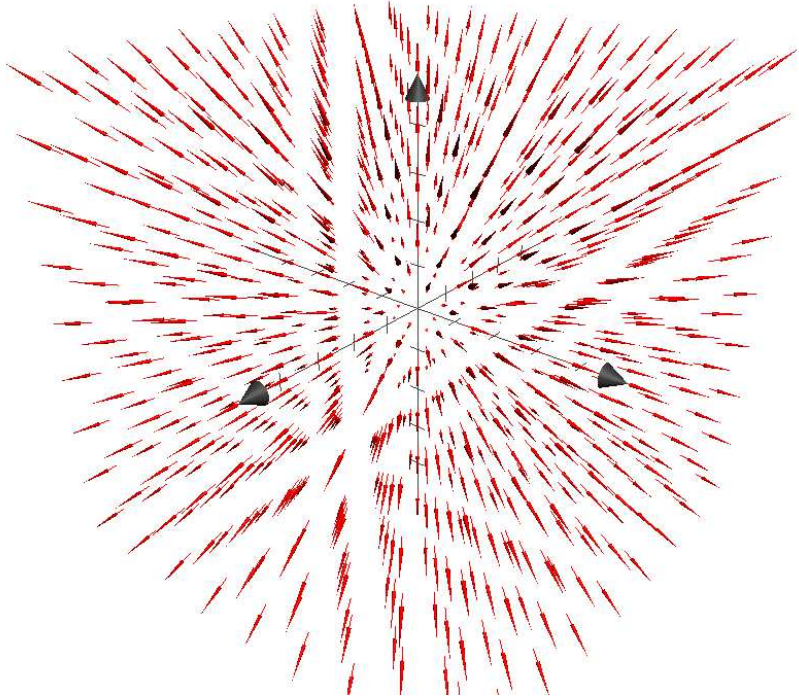} \end{center} \end{figure}

Let's see this result in action.  Consider the scalar function
\[ f(x,y,z) = \sqrt{x^2 + y^2 + z^2}. \]
\noindent We compute the gradient using the Chain Rule:
\begin{align*} \frac {\partial f}{\partial x} & = \frac 12 \left( x^2 + y^2 + z^2 \right)^{-1/2} \cdot 2 \, x = \frac {x}{\sqrt{x^2 + y^2 + z^2}}, \\ \frac {\partial f}{\partial y} & = \frac 12 \left( x^2 + y^2 + z^2 \right)^{-1/2} \cdot 2 \, y = \frac {y}{\sqrt{x^2 + y^2 + z^2}}, \\ \frac {\partial f}{\partial z} & = \frac 12 \left( x^2 + y^2 + z^2 \right)^{-1/2} \cdot 2 \, z = \frac {z}{\sqrt{x^2 + y^2 + z^2}}. \end{align*}
\noindent This gives the vector field
\[ \mathbf{F} = \nabla f = \frac {x \, \mathbf{i} + y \, \mathbf{j} + z \, \mathbf{k}}{\sqrt{x^2 + y^2 + z^2}}. \]
\noindent  (See Figure \ref{vector_field_1} for a plot.  To keep with the gastronomical motif of this article, perhaps we should call this direction field a ``Prickly Pear''?)  Recall that the curl of this vector field is
\[ \nabla \times \mathbf{F} = \left(\frac {\partial^2 f}{\partial y \, \partial z}  - \frac {\partial^2 f}{\partial z \, \partial y}\right)\mathbf{i}+\left(\frac {\partial^2 f}{\partial z \, \partial x} - \frac {\partial^2 f}{\partial x \, \partial z}\right)\mathbf{j}+\left(\frac {\partial^2 f}{\partial x \, \partial y} - \frac {\partial^2 f}{\partial y \, \partial x}\right)\mathbf{k}. \]
\noindent  We compute the mixed partial derivatives as follows:
\begin{align*}
\frac {\partial^2 f}{\partial y \, \partial z} & = \frac {\partial}{\partial y} \left[ \frac {z}{\sqrt{x^2 + y^2 + z^2}} \right] = - \frac 12 \, z \, \left( x^2 + y^2 + z^2 \right)^{-3/2} \cdot 2 \, y \\ & = - \frac {y \, z}{\bigl(x^2 + y^2 + z^2 \bigr)^{3/2}} \\[5pt]
\frac {\partial^2 f}{\partial z \, \partial y} & = \frac {\partial}{\partial z} \left[ \frac {y}{\sqrt{x^2 + y^2 + z^2}} \right] = - \frac 12 \, y \, \left( x^2 + y^2 + z^2 \right)^{-3/2} \cdot 2 \, z \\ & = - \frac {y \, z}{\bigl( x^2 + y^2 + z^2 \bigr)^{3/2}}\\
\end{align*}
\noindent Note that the other mixed partial derivatives give rise to a similar function.  Thus, $\nabla \times \mathbf{F} = \mathbf{0}$.

\subsection*{Relating Curls and Divergence}

Let's return to the second of our motivating questions.  We'll show the following:

\
\begin{center} \textit{$\mathbf{G} = \nabla \times \mathbf{F}$ is a curl if and only if the divergence $\nabla \cdot \mathbf{G} = 0$.} \end{center}
\

\noindent First assume that $\mathbf{G} = \nabla \times \mathbf{F}$ for some $\mathbf{F} = M \, \mathbf{i} + N \, \mathbf{j} + P \, \mathbf{k}$.  That is,
\[ S \, \mathbf{i} + T \, \mathbf{j} + U \, \mathbf{k} =\left( \frac {\partial P}{\partial y} - \frac {\partial N}{\partial z} \right) \mathbf{i} + \left( \frac {\partial M}{\partial z} - \frac {\partial P}{\partial x} \right) \mathbf{j} + \left( \frac {\partial N}{\partial x} - \frac {\partial M}{\partial y} \right) \mathbf{k}. \]
\noindent We compute the divergence as
\begin{align*}
\nabla \cdot \mathbf{G} & = \frac {\partial S}{\partial x} + \frac {\partial T}{\partial y} +\frac {\partial U}{\partial z}  \\
&=\left( \frac {\partial^2 P}{\partial x \, \partial y} - \frac {\partial^2 N}{\partial x \, \partial z} \right) + \left(\frac {\partial^2 M}{\partial y \, \partial z} - \frac {\partial^2 P}{\partial y \, \partial x}\right) +\left( \frac {\partial^2 N}{\partial z \, \partial x} - \frac {\partial^2 M}{\partial z \, \partial y} \right) \\
& = 0.
\end{align*}
\noindent (Here, we assume that $M$, $N$, and $P$ are continuously twice-differentiable so that we can interchange the order of differentiation.)\\

Now what if $\nabla \cdot \mathbf{G} = 0$?  In this case, we write
\[ \frac {\partial S}{\partial x} + \frac {\partial T}{\partial y} + \frac {\partial U}{\partial z} = 0 \qquad \implies \qquad U(x,y,z) = U(x,y,z_0) - \int_{z_0}^z \left[ \frac {\partial S}{\partial x} + \frac {\partial T}{\partial y} \right] d \zeta \]
\noindent for any fixed point $x_0 \, \mathbf{i} + y_0 \, \mathbf{j} + z_0 \, \mathbf{k} \in D$.   Define the function $\mathbf{F} = M \, \mathbf{i} + N \, \mathbf{j} + P \, \mathbf{k}$ in terms of the definite integrals
\[ \begin{aligned} M(x,y,z) & = \int_{z_0}^z T(x,y,\zeta) \, d \zeta - \int_{y_0}^y U(x, \tau, z_0) \, d \tau \\[5pt] N(x,y,z) & = - \int_{z_0}^z S(x,y,\zeta)  \, d \zeta \\[5pt] P(x,y,z) & = 0. \end{aligned} \]
\noindent Upon changing the order of differentiation and integration, we have the derivatives
\[ \begin{aligned} \frac {\partial P}{\partial y} - \frac {\partial N}{\partial z} & = S(x,y,z) \\[5pt] \frac {\partial M}{\partial z} - \frac {\partial P}{\partial x} & = T(x,y,z) \\[5pt]  \frac {\partial N}{\partial x} - \frac {\partial M}{\partial y} & = U(x,y,z_0) - \int_{z_0}^z \left[ \frac {\partial S}{\partial x} + \frac {\partial T}{\partial y}  \right] d \zeta = U(x,y,z). \end{aligned} \]
\noindent It follows that $\nabla \times \mathbf{F} = \mathbf{G}$.

\subsection*{Exercise}

\begin{figure}[t] \begin{center} \caption{Plot of $\mathbf{G} = x \, \mathbf{i} + 2 \, y \, \mathbf{j} -3 \, z \, \mathbf{k}$} \label{vector_field_2} \includegraphics[width=0.85\textwidth]{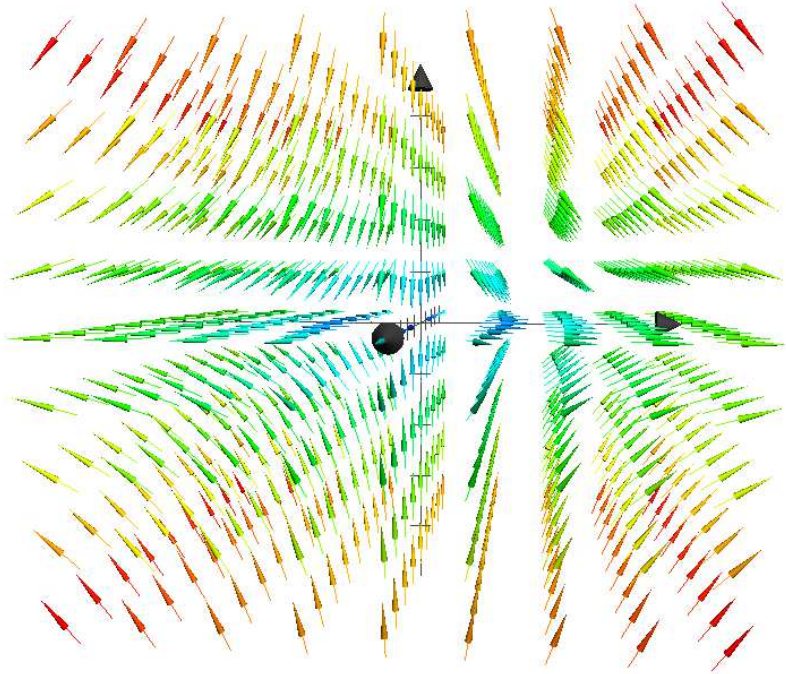} \end{center} \end{figure}

Fix three real numbers $a$, $b$, and $c$ such that $a + b + c = 0$, and define the function
\[ \mathbf{G} = a \, x \, \mathbf{i} + b \, y \, \mathbf{j} + c \, z \, \mathbf{k}. \]
\noindent (See Figure \ref{vector_field_2} for a plot.)  Show that $\nabla \cdot \mathbf{G} = 0$.  Moreover, find a function $\mathbf{F}$ such that $\mathbf{G} = \nabla \times \mathbf{F}$.  (\emph{Hint:} Guess a function in the form $\mathbf{F} = A \, y \, z \, \mathbf{i} + B \, x \, z \, \mathbf{j} + C \, x \, y \, \mathbf{k}$ for some real numbers $A$, $B$, and $C$.)  How many such functions $\mathbf{F}$ do you think there are?

\subsection*{Example}

Consider the vector function
\[ \mathbf{G} = \nabla f = \frac {x \, \mathbf{i} + y \, \mathbf{j} + z \, \mathbf{k}}{\sqrt{x^2 + y^2 + z^2}} \qquad \text{in terms of} \qquad f(x,y,z) = \sqrt{x^2 + y^2 + z^2}. \]

\noindent (Remember the ``Prickly Pear''?)  We will show that there is \emph{no} vector potential $\mathbf{F}$ such that $\mathbf{G} = \nabla \times \mathbf{F}$.  The idea is to show that the divergence of $\mathbf{G}$ is nonzero.  To this end, we compute the higher-order partial derivatives using the Quotient Rule:
\begin{align*}
\frac {\partial^2 f}{\partial x^2} & = \frac {\partial}{\partial x} \left[ \frac {x}{\sqrt{x^2 + y^2 + z^2}} \right] = \frac {y^2 + z^2}{\left( x^2 + y^2 + z^2 \right)^{3/2}} \\[5pt]
\frac {\partial^2 f}{\partial y^2} & = \frac {\partial}{\partial y} \left[ \frac {y}{\sqrt{x^2 + y^2 + z^2}} \right] = \frac {x^2 + z^2}{\left( x^2 + y^2 + z^2 \right)^{3/2}} \\[5pt]
\frac {\partial^2 f}{\partial z^2} & = \frac {\partial}{\partial z} \left[ \frac {z}{\sqrt{x^2 + y^2 + z^2}} \right] = \frac {x^2 + y^2}{\left( x^2 + y^2 + z^2 \right)^{3/2}} \\[5pt]
\end{align*}
\noindent (Check as an exercise!)  Hence we find the divergence
\[ \nabla \cdot \mathbf{G} = \frac {\partial^2 f}{\partial x^2} + \frac {\partial^2 f}{\partial y^2} + \frac {\partial^2 f}{\partial z^2} = \frac {2 \, x^2 + 2 \, y^2 + 2 \, z^2}{\left( x^2 + y^2 + z^2 \right)^{3/2}} = \frac {2}{f}. \]
\noindent Since $\nabla \cdot \mathbf{G} \neq 0$, there cannot exist a function $\mathbf{F}$ such that $\mathbf{G} = \nabla \times \mathbf{F}$.

\section*{Differential Forms}

Naturally, the three motivating questions in the introduction are rather naive questions to ask  -- although a bit difficult to answer! -- so let's make them a little more interesting.  We'll rephrase the definitions above in a way using differential $k$-forms; this will make integration more natural.  \\

\begin{itemize}
\item A \emph{0-form} is a continuously differentiable function $f: D \to \mathbb R$.  These are the functions we know and love from the one-variable Differential Calculus.\\
\item A \emph{1-form} is an arc length differential in the form $\eta = M \, dx + N \, dy + P \, dz$ for some continuously differentiable vector field $\mathbf{F}: D \to \mathbb R^3$ in the form $\mathbf{F} = M \, \mathbf{i} + N \, \mathbf{j} + P \, \mathbf{k}$.  Note that we can express this in the form $\eta = \mathbf{F} \cdot d \mathbf{r}$ using the dot product.   We'll use this to define path integrals later.  \\
\item A \emph{2-form} is an area differential in the form $\omega = S \, dy \, dz + T \, dx \, dz + U \, dx \, dy$ for some continuously differentiable vector field $\mathbf{G}: D \to \mathbb R^3$ in the form $\mathbf{G} = S \, \mathbf{i} + T \, \mathbf{j}  + U \, \mathbf{k}$.  Note that we can express this in the form $\omega = \mathbf{G} \cdot d \mathbf{A}$.  As the notation suggests, we'll use this to define surface integrals.  \\
\item \emph{3-form} is a volume differential in the form $\nu = f \, dx \, dy \, dz$ for some continuously differentiable function $f: D \to \mathbb R$.  Note that we can express this in the form $\nu = f \, dV$.  As the notation suggests, we'll use this to define volume integrals.
\end{itemize}

We'll denote $\Omega^k(D)$ as the collection of $k$-forms on $D$ for $k = 0$, 1, 2, and 3.    Boothby \cite{MR861409} provides a rigorous treatment of differential forms. \\

We claim that each $\Omega^k(D)$ these is a linear vector space.  This means that the linear combination of two $k$-forms is another $k$-form.  Consider, for example, 1-forms and 2-forms.  Given scalars $\alpha$ and $\beta$, we have the identities
\[ \left. \begin{aligned} \alpha \, & \bigl( M_1 \, dx + N_1 \, dy + P_1 \, dz \bigr) + \beta \, \bigl( M_2 \, dx + N_2 \, dy + P_2 \, dz \bigr) \\ & \qquad \qquad = \bigl( \alpha \, M_1 + \beta \, M_2 \bigr) \, dx + \bigl( \alpha \, N_1 + \beta \, N_2 \bigr) \, dy + \bigl( \alpha \, P_1 + \beta \, P_2 \bigr) \, dy  \\[5pt]  \alpha \, & \bigl( S_1 \, dy \, dz + T_1 \, dx \, dz + U_1 \, dx  \, dy \bigr) + \beta \, \bigl( S_2 \, dy \, dz + T_2 \, dx \, dz + U_2 \, dx  \, dy \bigr) \\ & \qquad \qquad = \bigl( \alpha \, S_1 + \beta \, S_2 \bigr) \, dy \, dz + \bigl( \alpha \, T_1 + \beta \, T_2 \bigr) \, dx \, dz + \bigl( \alpha \, U_1 + \beta \, U_2 \bigr) \, dx \, dy \end{aligned} \right.\]

\noindent Thus $\Omega^k(D)$ is a linear vector space for $k=1, 2$.  The argument for $k = 0, 3$ is similar.

\subsection*{Maps Between Differential Forms}

How are these linear vector spaces related?  Well, there are linear maps between them!  Consider the following ``differential'' map $d: \Omega^0(D) \to \Omega^1(D)$ from 0-forms to 1-forms:
\[ \begin{matrix} f \\[5pt] \downarrow \\[5pt] df =  \dfrac {\partial f}{\partial x} \, dx + \dfrac {\partial f}{\partial y} \, dy + \dfrac {\partial f}{\partial z} \, dz = \nabla f \cdot d \mathbf{r}. \end{matrix} \]
\noindent This is linear because $d \bigl( \alpha \, f_1 + \beta \, f_2 \bigr) = \alpha \, d f_1 + \beta \, d f_2$.  Here, we write $d \mathbf{r} = dx \, \mathbf{i} + dy \, \mathbf{j} + dz \, \mathbf{k}$ as the arc length differential -- expressed as a vector.  Hence we find the answer to our first motivating question:

\
\begin{center} \textit{The vector-valued function $\mathbf{F} = \nabla f$ has a scalar potential $f$ if and only if the 1-form $\mathbf{F} \cdot d \mathbf{r} = df$ is the differential of some 0-form $f$.} \end{center}
\

Similarly, consider the following ``differential'' map  $d: \Omega^1(D) \to \Omega^2(D)$ from 1-forms to 2-forms:
\[ \begin{matrix} \eta = M \, dx + N \, dy + P \, dz = \mathbf{F} \cdot d \mathbf{r} \\[5pt] \downarrow \\[5pt] \begin{aligned} d \eta & =  \left( \dfrac {\partial P}{\partial y} - \dfrac {\partial N}{\partial z} \right) dy \, dz + \left( \dfrac {\partial M}{\partial z} - \dfrac {\partial P}{\partial x} \right) dx \, dz + \left( \dfrac {\partial N}{\partial x} - \dfrac {\partial M}{\partial y} \right) dx \, dy \\ & = \bigl( \nabla \times \mathbf{F} \bigr) \cdot d \mathbf{A}. \end{aligned} \end{matrix} \]
\noindent This is linear because $d \bigl( \alpha \, \eta_1 + \beta \, \eta_2 \bigr) = \alpha \, d \eta_1 + \beta \, d \eta_2$.   Here, we write $d \mathbf{A} = dy \, dz \, \mathbf{i} + dx \, dz \, \mathbf{j} + dx \, dy \, \mathbf{k}$ as the area differential -- also expressed as a vector.  Hence we find the answer to our second motivating question:

\
\begin{center} \textit{The vector-valued function $\mathbf{G} = \nabla \times \mathbf{F}$ has a vector potential $\mathbf{F}$  if and only if the 2-form $\mathbf{G} \cdot d \mathbf{A} = d \eta$ is the differential of some 1-form $\eta =  \mathbf{F} \cdot d \mathbf{r}$.} \end{center}
\

\noindent Recall that $\mathbf{F} = \nabla f$ is a gradient if and only if the curl $\nabla \times \mathbf{F} = \mathbf{0}$.  In other words, given a 1-form $\eta = M \, dx + N \, dy + P \, dz = \mathbf{F} \cdot d \mathbf{r}$, we have $\eta = d f$ as the differential of a 0-form if and only if the 2-form $d \eta = 0$.  This is yet another answer to our first question!

Finally, consider the following ``differential'' map  $d: \Omega^2(D) \to \Omega^3(D)$ from 2-forms to 3-forms:
\[ \begin{matrix} \omega = S \, dy \, dz + T \, dx \, dz + U \, dx \, dy = \mathbf{G} \cdot d \mathbf{A} \\[5pt] \downarrow \\[5pt] d \omega = \left( \dfrac {\partial S}{\partial x} + \dfrac {\partial T}{\partial y} + \dfrac {\partial U}{\partial z} \right) dx \, dy \, dz = \bigl( \nabla \cdot \mathbf{G} \bigr) \, dV. \end{matrix} \]
\noindent Again, this is linear because $d \bigl( \alpha \, \omega_1 + \beta \, \omega_2 \bigr) = \alpha \, d \omega_1 + \beta \, d \omega_2$.  Remember that $dV = dx \, dy \, dz$ is the volume differential.  Hence we find the answer to the last of our three motivating questions:

\
\begin{center} \textit{The vector-valued function $f = \nabla \cdot \mathbf{G}$ has a vector potential $\mathbf{G}$ if and only if the  3-form $f \, dV = d \omega$ is the differential of some 2-form $\omega =  \mathbf{G} \cdot d \mathbf{A}$.} \end{center}
\

\noindent Recall that $\mathbf{G} = \nabla \times \mathbf{F}$ is a curl if and only if the divergence $\nabla \cdot \mathbf{G} = 0$.  In other words, given a 2-form $\omega = S \, dy \, dz + T \, dx \, dz + U \, dx \, dy = \mathbf{G} \cdot d \mathbf{A}$,  we have $\omega = d \eta$ as the differential of a 1-form if and only if the 3-form $d \omega = 0$.  This is yet another answer to our second question! \\

We'll summarize this via the following diagram of maps:
\[ \begin{CD} \{ 0 \} @>>> \Omega^0(D) @>{\text{gradient}}>> \Omega^1(D) @>{\text{curl}}>> \Omega^2(D) @>{\text{divergence}}>> \Omega^3(D) @>>> \{ 0 \}. \end{CD} \]
\noindent All of the various definitions of differentiation for functions in three variables are all related via linear transformations on the vector spaces of $k$-forms!

\subsection*{Example}

\begin{figure}[h] \begin{center} \caption{Plot of $\mathbf{F} = y \, z \, \mathbf{i} + x \, z \, \mathbf{j} + x \, y \, \mathbf{k}$} \label{vector_field_3} \includegraphics[width=0.85\textwidth]{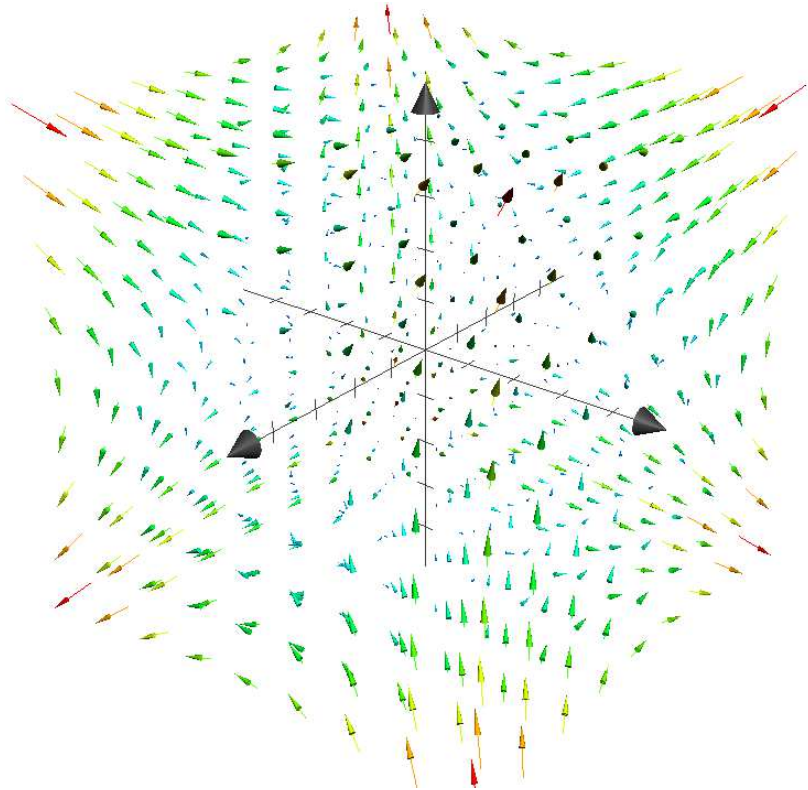} \end{center} \end{figure}

Let's consider the differential 1-form
\[ \eta = y \, z \, dx + x \, z \, dy + x \, y \, dz. \]
\noindent This is in the form $M \, dx + N \, dy + P \, dz$ for $M = y \, z$, $N = x \, z$, and $P = x \, y$.  (See Figure \ref{vector_field_3} for a plot of $\mathbf{F} = M \, \mathbf{i} + N \, \mathbf{j} + P \, \mathbf{k}$.)  We'll compute the differential $d \eta$.  We have the partial derivatives
\[  \dfrac {\partial P}{\partial y} - \dfrac {\partial N}{\partial z} = x - x = 0, \quad \dfrac {\partial M}{\partial z} - \dfrac {\partial P}{\partial x} = y - y = 0, \quad \dfrac {\partial N}{\partial x} - \dfrac {\partial M}{\partial y} = z - z = 0.  \]
\noindent Hence $d \eta = 0$.

\subsection*{Exercise}

Fix three real numbers $a$, $b$, and $c$ such that $a + b + c = 0$, and differential 2-form
\[ \omega = a \, x \, dy \, dz + b \, y \, dx \, dz + c \, z \, dx \, dy. \]
\noindent Show that $\omega = d \eta$ is the differential of the 1-form $\eta = b \, y \, z \, dx - a \, x \, z \, dy$.  Is this the only 1-form $\eta$ such that $\omega = d \eta$?

\section*{Integration}

We know now how the various forms of differentiation are all related, but what about integration?  First, we review what we mean by integration.  As $D \subseteq \mathbb R^3$ is a subset of three-dimensional space, we can define integration for one variable, two variables, or even three variables.

\subsection*{Fundamental Theorem of Calculus}

In order to define an integral for one variable, let $\mathbf{r}: I \to D$ be a continuously differentiable map, which we write in the form $\mathbf{r} = x \, \mathbf{i} + y \, \mathbf{j} + z \, \mathbf{k}$, defined on a closed interval $I = [a,b]$.  We say that the image $C \subseteq D$ is a \emph{path}.  Given vector field $\mathbf{F}: D \to \mathbb R^3$, the composition gives the 1-form $\eta = \mathbf{F}\bigl( \mathbf{r}(t) \bigr) \cdot d \mathbf{r} $, so naturally a \emph{path integral} as
\[ \int_{C} \eta = \int_a^b \left[ M \bigl( \mathbf{r}(t) \bigr) \, \frac {dx}{dt} + N \bigl( \mathbf{r}(t) \bigr) \, \frac {dy}{dt} + P \bigl( \mathbf{r}(t) \bigr) \, \frac {dz}{dt} \right] dt. \]
\noindent Recall that $\mathbf{F} = \nabla f$ has a scalar potential $f$ if and only if $\eta = d f$ is the differential of a 0-form $f$.  In this case, the path integral simplifies to
\[ \int_a^b \frac d{dt} \biggl[ f \bigl( \mathbf{r}(t) \bigr) \biggr] dt  = \int_{C} d f = f(\mathbf{r}) \biggr|_{\partial C} = f \bigl( \mathbf{r}(b) \bigr) - f \bigl( \mathbf{r}(a) \bigr). \]
\noindent Hence the integral is independent of the path $C$, and it is only dependent on the \emph{endpoints} $\partial C = \left \{ \mathbf{r}(a), \ \mathbf{r}(b) \right \}$.  This is just the \emph{Fundamental Theorem of Calculus}.

\subsection*{Example}

Let $D = \mathbb R^3$ be all of three-dimensional space, and let $C$ denote the circle of radius $r$ in the plane.  Note that $C$ is the set of points $x \, \mathbf{i} + y \, \mathbf{j} + z \, \mathbf{k}$ such that $x^2 + y^2 = r^2$ and $z = 0$.  We may think of this as being the image of the map $\mathbf{r}: [0, 2 \pi] \to \mathbb R^3$ which sends $t \mapsto r \, \cos t \, \mathbf{i} + r \, \sin t \, \mathbf{j}$.  The circle is just an example of a path!

\begin{figure}[h] \begin{center} \caption{Plot of $\mathbf{F} = - (y/2) \, \mathbf{i} + (x/2) \, \mathbf{j}$} \label{vector_field_4} \includegraphics[width=0.85\textwidth]{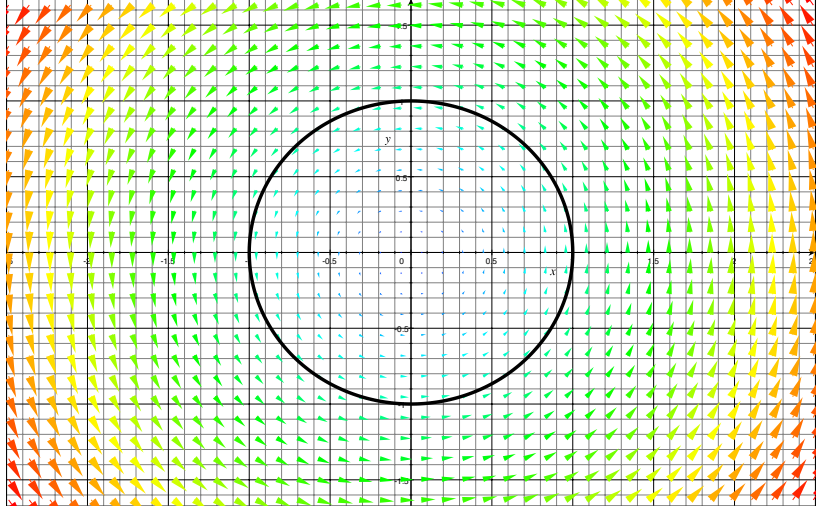} \end{center} \end{figure}

Now let's consider a vector field $\mathbf{F} = - (y/2) \, \mathbf{i} + (x/2) \, \mathbf{j}$.  (See Figure \ref{vector_field_4} for a plot.)  We'll explain why there does \emph{not} exist a scalar function $f$ such that $\mathbf{F} = \nabla f$.  The idea is to suppose that $f$ does indeed exist and then compute the integral of $\mathbf{F}$ around the path $C$.  Let's consider the following differential 1-form:
\begin{align*} \eta & = \mathbf{F}\bigl( \mathbf{r}(t) \bigr) \cdot d \mathbf{r} = \biggl( - \frac {r \, \cos t}2 \, \mathbf{i} + \frac {r \, \sin t}2 \, \mathbf{j} \biggr) \cdot \biggl( -r \, \cos t \, dt \, \mathbf{i} + r \, \sin t \, dt \, \mathbf{j} \biggr) \\ & = \frac {r^2 \, \cos^2 t \, dt + r^2 \, \sin^2 t \, dt}2 = \frac {r^2}2 \, dt. \end{align*}
\noindent If $\mathbf{F} = \nabla f$, then $\int_C \eta = f \bigl( \mathbf{r}(2\pi) \bigr) - f \bigl( \mathbf{r}(0) \bigr) = 0.$  But actually, 
\[ \int_C \eta = \int_0^{2 \pi} \frac {r^2}{2} \, dt = \pi \, r^2. \]
\noindent Hence there is no function $f$ such that $\mathbf{F} = \nabla f$.  (Of course, we could have seen this sooner by computing the curl $\nabla \times \mathbf{F} = \mathbf{k}$ and realizing it's a nonzero vector.)

\subsection*{Stokes' Theorem and Green's Theorem}

In order to define an integral for two variables, let $\mathbf{r}: R \to D$ be a continuously differentiable map, which we write in the form $\mathbf{r} = x \, \mathbf{i} + y \, \mathbf{j} + z \, \mathbf{k}$, that is defined on an closed region $R = [a,b] \times [c,d]$.  We say that the image $S \subseteq D$ is a \emph{surface}, and assume that its boundary $C = \partial S$ is a curve as before.  Given vector field $\mathbf{G}: D \to \mathbb R^3$, the composition gives the 2-form $\omega = \mathbf{G} \bigl( \mathbf{r}(u,v) \bigr) \cdot d \mathbf{A}$, so naturally \emph{surface integral} is defined as
\[ \iint_{S} \omega = \int_a^b \int_c^d \left| \begin{matrix} S \bigl( \mathbf{r}(u,v) \bigr) & T \bigl( \mathbf{r}(u,v) \bigr) & U \bigl( \mathbf{r}(u,v) \bigr) \\[8pt] \dfrac {\partial x}{\partial u} & \dfrac {\partial y}{\partial u} & \dfrac {\partial z}{\partial u} \\[8pt] \dfrac {\partial x}{\partial v} & \dfrac {\partial y}{\partial v} & \dfrac {\partial z}{\partial v} \end{matrix} \right| du \, dv. \]
\noindent Recall that $\mathbf{G} = \nabla \times \mathbf{F}$ has a vector potential $\mathbf{F}$ if and only if $\omega = d \eta$ is the differential of a 1-form $\eta$.  In this case, the surface integral simplifies to
\[  \iint_{S} \bigl( \nabla \times \mathbf{F} \bigr) \cdot d \mathbf{A} = \iint_{S} d \eta = \int_{\partial S} \eta = \int_{C} \mathbf{F} \cdot d \mathbf{r}. \]
\noindent Hence the integral is independent of the surface $\mathcal{S}$, and it is only dependent on the boundary $C = \partial S$.  This is known as \emph{Stokes' Theorem}.

Let's consider a special case, where $\mathbf{r}: R \to \mathbb R^2$ actually maps into the plane.  Then $S  \subseteq \mathbb R^2$, so write $\mathbf{F} = M \, \mathbf{i} + N \, \mathbf{j}$.  Stokes' Theorem in the plane reduces to the statement
\[ \iint_{S} \left( \frac {\partial N}{\partial x} - \frac {\partial M}{\partial y} \right) dx \, dy = \iint_{S} d \eta = \int_{\partial R} \eta = \int_{\partial S} \bigl[ M \, dx + N \, dy \bigr]. \]
\noindent This is known as \emph{Green's Theorem}.  Similarly, for the \emph{orthogonal} vector field $\mathbf{F}_{\perp} = -N \, \mathbf{i} + M \, \mathbf{j}$, the \emph{Divergence Theorem in the Plane} is the expression
\[ \iint_{S} \left( \frac {\partial M}{\partial x} + \frac {\partial N}{\partial y} \right) dx \, dy = \int_{\partial S} \bigl[ -N \, dx + M \, dy \bigr]. \]

\subsection*{Example}

Let $D = \mathbb R^3$ be all of three-dimensional space, and let $S$ denote the disk of radius $r$ in the plane.  The latter is just the set of points $x \, \mathbf{i} + y \, \mathbf{j} + z \, \mathbf{k}$ such that $x^2 + y^2 \leq r^2$ and $z = 0$.  We may think of this as being the image of the map $\mathbf{r}: [0,r] \times [0, 2 \pi] \to \mathbb R^3$ which sends $(\rho,\theta) \mapsto \rho \, \cos \theta \, \mathbf{i} + \rho \, \sin \theta \, \mathbf{j}$.  Note that the boundary $C = \partial S$ is simply the circle $x^2 + y^2 = r^2$, which we considered before.

Let's return to the vector field $\mathbf{F} = - (y/2) \, \mathbf{i} + (x/2) \, \mathbf{j}$.   We'll consider the differential 2-form $\omega = \bigl( \nabla \times \mathbf{F} \bigr) \cdot d \mathbf{A} = \mathbf{k} \cdot d \mathbf{A} = dx \, dy$.  Stokes' Theorem (or really Green's Theorem, since it's in the plane) states that
\[ \iint_S \omega = \int_C \mathbf{F} \cdot d \mathbf{r} = \pi \, r^2. \]

\noindent (Recall the previous example.)  Let's compute the integral on the left-hand side in a different way.  This differential 2-form can also expressed in \emph{polar coordinates} using the determinant
\begin{align*} \omega & = \left| \begin{matrix} S \bigl( \mathbf{r}(\rho,\theta) \bigr) & T \bigl( \mathbf{r}(\rho,\theta) \bigr) & U \bigl( \mathbf{r}(\rho,\theta) \bigr) \\[8pt] \dfrac {\partial x}{\partial \rho} & \dfrac {\partial y}{\partial \rho} & \dfrac {\partial z}{\partial \rho} \\[8pt] \dfrac {\partial x}{\partial \theta} & \dfrac {\partial y}{\partial \theta} & \dfrac {\partial z}{\partial \theta} \end{matrix} \right| d \rho \, d \theta = \left| \begin{matrix} 0 & 0 & 1 \\ \cos \theta & \sin \theta & 0 \\ - \rho \, \sin \theta & \rho \, \cos \theta & 0 \end{matrix} \right| d \rho \, d \theta \\[8pt] & = \rho \, d \rho \, d \theta. \end{align*}

\noindent Hence we have the integral
\[ \iint_S dx \, dy = \iint_S \omega = \int_0^r \int_0^{2 \pi} \rho \, d \rho \, d \theta = \left[ \int_0^r  \rho \, d \rho \right] \left[ \int_0^{2 \pi} d \theta \right] = \pi \, r^2. \]
\noindent Of course, this is just the area of the disk $S$.

\subsection*{Exercise}

Let $S$ be any surface in the plane $\mathbb R^2$ with a boundary $C = \partial S$.  Show that its area can be computed using a path integral around along $C$.  That is,
\[ \text{Area}(S) = \iint_S dx \, dy = \int_{C} \frac {-y \, dx + x \, dy}{2}.\]

\subsection*{Gauss' Theorem and the Divergence Theorem.}

In order to define an integral for three variables, let $\mathbf{r}: B \to \mathbb R^3$ be a continuously differentiable map, which we write in the form $\mathbf{r} = x \, \mathbf{i} + y \, \mathbf{j} + z \, \mathbf{k}$, defined on an closed region $B = [a,b] \times [c,d] \times [p,q]$.   We say that the image $D \subseteq \mathbb R^3$ is a \emph{region}, and assume that its boundary $S = \partial D$ is a surface as before.  Given a scalar function $f: D \to \mathbb R$, the composition gives the differential 3-form $\nu = f\bigl( \mathbf{r}(u,v,w) \bigr) \, dV$, so naturally define a \emph{volume integral} as
\[ \iiint_{D} \nu = \int_a^b \int_c^d \int_p^q f \bigl( \mathbf{r}(u,v,w) \bigr) \left| \begin{matrix} \dfrac {\partial x}{\partial u} & \dfrac {\partial y}{\partial u} & \dfrac {\partial z}{\partial u} \\[8pt] \dfrac {\partial x}{\partial v} & \dfrac {\partial y}{\partial v} & \dfrac {\partial z}{\partial v} \\[8pt] \dfrac {\partial x}{\partial w} & \dfrac {\partial y}{\partial w} & \dfrac {\partial z}{\partial w} \end{matrix} \right| du \, dv \, dw. \]

\noindent (As with surface integrals, we've expressed the integrand above using a $3 \times 3$ determinant.  Most texts refer to this as the \emph{Jacobian} of the transformation $\mathbf{r}: B \to D$.)  Recall that $f = \nabla \cdot \mathbf{G}$ has a vector potential $\mathbf{G}$ if and only if $\nu = d \omega$ is a differential of a 2-form $\omega$.  In this case, the volume integral simplifies to
\[  \iiint_{D} \bigl( \nabla \cdot \mathbf{G} \bigr) \, dV = \iiint_{D} d \omega = \iint_{\partial D} \omega = \iint_{S} \mathbf{G} \cdot d \mathbf{A}. \]
\noindent Hence the integral is independent of the region $V$, and it is only dependent on the boundary $S = \partial D$.  This is known as \emph{Gauss' Theorem}, or the \emph{Divergence Theorem}.

\subsection*{Example}

Let $D$ denote the solid sphere of radius $r$ in three-dimensional space.  This is just the set of points $x \, \mathbf{i} + y \, \mathbf{j} + z \, \mathbf{k}$ such that $x^2 + y^2 + z^2 \leq r^2$.  We may think of this as being the image of the map $\mathbf{r}:  [0,r] \times [0, \pi]  \times [0, 2 \pi] \to \mathbb R^3$ which sends
\[ (\rho, \phi,\theta) \mapsto \rho \, \sin \phi \, \cos \theta \, \mathbf{i} + \rho \, \sin \phi \, \sin \theta \, \mathbf{j} + \rho \, \cos \phi \, \mathbf{k}. \]
\noindent Note that the boundary $S = \partial D$ is simply the sphere $x^2 + y^2 + z^2 = r^2$.

Let's return to the vector field
\[ \mathbf{G} = \nabla g = \frac {x \, \mathbf{i} + y \, \mathbf{j} + z \, \mathbf{k}}{\sqrt{x^2 + y^2 + z^2}} \qquad \text{in terms of} \qquad g(x,y,z) = \sqrt{x^2 + y^2 + z^2}. \]

\noindent (The ``Prickly Pear'' again!)  We will compute integrals and verify the Divergence Theorem.  We may express the volume element $dV = dx \, dy \, dz$ in \emph{spherical coordinates} using the determinant
\begin{align*} dV & = \left| \begin{matrix} \dfrac {\partial x}{\partial \rho} & \dfrac {\partial y}{\partial \rho} & \dfrac {\partial z}{\partial \rho}  \\[8pt] \dfrac {\partial x}{\partial \phi} & \dfrac {\partial y}{\partial \phi} & \dfrac {\partial z}{\partial \phi} \\[8pt] \dfrac {\partial x}{\partial \theta} & \dfrac {\partial y}{\partial \theta} & \dfrac {\partial z}{\partial \theta} \end{matrix} \right| d\rho \, d \phi \, d \theta = \left| \begin{matrix} \cos \theta \, \sin \phi & \sin \theta \, \sin \phi & \cos \phi  \\[5pt] \rho \, \cos \theta \, \cos \phi & \rho \, \sin \theta \, \cos \phi & - \rho \, \sin \phi \\[5pt] - \rho \, \sin \theta \, \sin \phi & \rho \, \cos \theta \, \sin \phi & 0 \end{matrix} \right| d\rho \, d \phi \, d \theta \\[5pt] & = \rho^2 \, \sin \phi \, d \rho \, d \phi \, d \theta. \end{align*}
\noindent Hence we have the differential 3-form
\[ \nu = \nabla \cdot \mathbf{G} \bigr( \mathbf{r}(\rho, \phi, \theta) \bigl) \, dV = \frac {2}{\rho} \cdot \bigl( \rho^2 \, \sin \phi \, d \rho \, d \phi \, d \theta \bigr) = 2 \, \rho \, \sin \phi \, d \rho \, d \phi \, d \theta. \]
\noindent (Recall that $\nabla \cdot \mathbf{G} = 2/g$.)  This gives the integral
\begin{align*} \iiint_{D} \nu & = \int_0^r \int_0^{\pi} \int_0^{2 \pi} 2 \, \rho \, \sin \phi \, d \rho \, d \phi \, d \theta = \left[ \int_0^r 2 \, \rho\, d \rho \right] \left[ \int_0^{\pi} \sin \phi \, d \phi \right] \left[ \int_0^{2 \pi}  d \theta \right] \\[5pt] & = 4 \, \pi \, r^2. \end{align*}

\noindent On the other hand, we may express the area element is the determinant
\begin{align*} d \mathbf{A} & = \left| \begin{matrix} \mathbf{i} & \mathbf{j} & \mathbf{k} \\[8pt] \dfrac {\partial x}{\partial \phi} & \dfrac {\partial y}{\partial \phi} & \dfrac {\partial z}{\partial \phi} \\[8pt] \dfrac {\partial x}{\partial \theta} & \dfrac {\partial y}{\partial \theta} & \dfrac {\partial z}{\partial \theta} \end{matrix} \right| d \phi \, d \theta = \left| \begin{matrix} \mathbf{i} & \mathbf{j} & \mathbf{k}  \\[5pt] \rho \, \cos \theta \, \cos \phi & \rho \, \sin \theta \, \cos \phi & - \rho \, \sin \phi \\[5pt] - \rho \, \sin \theta \, \sin \phi & \rho \, \cos \theta \, \sin \phi & 0 \end{matrix} \right| d \phi \, d \theta \\[5pt] & = \rho \, \sin \phi \, \bigl( \rho \, \sin \phi \, \cos \theta \, \mathbf{i} + \rho \, \sin \phi \, \sin \theta \, \mathbf{j} + \rho \, \cos \phi \, \mathbf{k} \bigr) \, d \phi \, d \theta. \end{align*}
\noindent Hence we have the differential 2-form
\[ \omega = \mathbf{G} \bigl( \mathbf{r}(\rho, \phi, \theta) \bigr) \cdot d \mathbf{A} = \frac {\mathbf{r}}{\rho} \cdot \bigl( \rho \, \sin \phi \, \mathbf{r} \, d \theta \, d \phi \bigr) = \rho^2 \, \sin \phi \, d \phi \, d \theta. \]
\noindent This gives the integral
\[ \iint_{\partial D} \omega = \int_0^{\pi} \int_0^{2 \pi} r^2 \, \sin \phi \, d \theta \, d \phi = r^2 \left[ \int_0^{\pi} \sin \phi \, d \phi \right] \left[ \int_0^{2 \pi} d \theta \right] = 4 \, \pi \, r^2. \]

\noindent (Remember that $\rho = r$ along the boundary $S = \partial D$.)  This indeed verifies that
\[  \iiint_{D} \bigl( \nabla \cdot \mathbf{G} \bigr) \, dV = \iint_{S} \mathbf{G} \cdot d \mathbf{A}. \]
\noindent Arfken and Weber \cite{MR1810939} provides a plethora of formulas for integrating in different coordinate systems other than spherical.

\subsection*{Exercise}

Let $D$ be any region in $\mathbb R^3$ with a boundary $S = \partial D$.  Show that its volume can be computed using a surface integral on $S$.  That is,
\[ \text{Vol}(D) = \iiint_D dx \, dy \, dz = \iint_{S} \frac {x \, dy \, dz + y \, dx \, dz + z \, dx \, dy}{3}.\]

\section*{The Moral of the Story}

In this paper we have shown that all of the theorems -- those for differentiation and integration -- can be expressed using differential forms.   Here's the idea in a nutshell for \emph{any} region $D$.  We consider a series of linear vector spaces $\Omega^k(D)$, where we have ``differential'' maps
\[ \begin{CD} \cdots @>>> \Omega^{k-1}(D) @>{d}>> \Omega^k(D) @>{d}>> \Omega^{k+1}(D) @>>> \cdots \end{CD} \]
\noindent  We would like to know the answer to the following question involving differentiation:
\vspace{-.10in}
\begin{center} \textit{Which $k$-forms $\omega \in \Omega^{k+1}(D)$ are in the form \\ $\omega = d \eta$ for some $(k-1)$-form $\eta \in \Omega^k(D)$?} \end{center}
\vspace{.05in}
\noindent The \emph{partial} answer should be:
\vspace{.05in}
\begin{center} \textit{If $\omega = d \eta$ for some $(k-1)$-form $\eta$, \\ then $d \omega = 0$ as a $(k+1)$-form.} \end{center}
\vspace{.05in}
\noindent A complete answer involves computing something called \emph{de Rham Cohomology}.  (The desert after the Main Course?)  For instance, we have computed de Rham Cohomlogy in this article for certain subsets of three-dimensional space.  We expect to generalize the integration formulas above by saying something like
\[ \int_{D} \omega = \int_{D} d \eta = \int_{\partial D} \eta \]
\noindent so that the integral would be independent of the region $D$, and would be only dependent on the boundary $\partial D$.  This is the \emph{Generalized Stokes' Theorem}.  Surprisingly, this entire theory can be worked out for ``many'' sets $D$.  But don't take our words for it -- just a friendly Differential Geometer!

\section*{Appendix: $3 \times 3$ Determinants}

A \emph{$3 \times 3$ determinant} is the quantity
\[ \left| \begin{tabular}{ccc} $a_{11}$ & $a_{12}$ & $a_{13}$ \\ $a_{21}$ & $a_{22}$ & $a_{23}$ \\ $a_{31}$ & $a_{32}$ & $a_{33}$ \end{tabular} \right| = \begin{aligned} a_{11}  (\, a_{22} \, a_{33} & - a_{23} \, a_{32}) - a_{12}( \, a_{23} \, a_{31} - \, a_{21} \, a_{33}) \\ & + a_{13}( \, a_{21} \, a_{32} - a_{22} \, a_{31}). \end{aligned} \]

\noindent Another way is to compute the determinant is by using the following diagram:
\[ \centerline{\xymatrix@!0{ \ar@{->}'[dr]'[ddrr]'[dddrrr] [ddddrrrr] & \ar@{->}'[dr]'[ddrr]'[dddrrr] [ddddrrrr] & \ar@{->}'[dr]'[ddrr]'[dddrrr] [ddddrrrr] & & \ar@{->}'[dl]'[ddll]'[dddlll] [ddddllll]  & \ar@{->}'[dl]'[ddll]'[dddlll] [ddddllll]  & \ar@{->}'[dl]'[ddll]'[dddlll] [ddddllll] \\ & a_{11} & a_{12} & a_{13} & a_{11} & a_{12} & \\  & a_{21} & a_{22} & a_{23} & a_{21} & a_{22} & \\ & a_{31} & a_{32} & a_{33} & a_{31} & a_{32} & \\ & & & & & &\\ }} \]
\noindent Here we multiply along the arrows, then add or subtract depending on the direction of the arrow:
\[ \left| \begin{tabular}{ccc} $a_{11}$ & $a_{12}$ & $a_{13}$ \\ $a_{21}$ & $a_{22}$ & $a_{23}$ \\ $a_{31}$ & $a_{32}$ & $a_{33}$ \end{tabular} \right| = \begin{aligned} & \left( a_{11} \, a_{22} \, a_{33} + a_{12} \, a_{23} \, a_{31} + a_{13} \, a_{21} \, a_{32} \right) \\ & \quad \qquad \qquad - \left( a_{11} \, a_{23} \, a_{32} +  a_{12} \, a_{21} \, a_{33} + a_{13} \, a_{22} \, a_{31} \right). \end{aligned} \]


\begin{thebibliography}{1}

\bibitem{MR1810939}
George~B. Arfken and Hans~J. Weber.
\newblock {\em Mathematical methods for physicists}.
\newblock Harcourt/Academic Press, Burlington, MA, fifth edition, 2001.

\bibitem{MR861409}
William~M. Boothby.
\newblock {\em An introduction to differentiable manifolds and {R}iemannian
  geometry}, volume 120 of {\em Pure and Applied Mathematics}.
\newblock Academic Press Inc., Orlando, FL, second edition, 1986.

\bibitem{Hughes-Hallett2004Calculus:-Singl}
Deborah Hughes-Hallett, Andrew~M. Gleason, William~G. McCallum, Daniel~E.
  Flath, Patti~Frazer Lock, Thomas~W. Tucker, David~O. Lomen, David Lovelock,
  David Mumford, Brad~G. Osgood, Douglas Quinney, Karen Rhea, and Jeff
  Tecosky-Feldman.
\newblock {\em Calculus: Single and Multivariable}.
\newblock Wiley-Interscience [John Wiley \& Sons], 4th edition, December 2004.

\bibitem{Larson2007Calculus:-Early}
Ron Larson, Robert~P. Hostetler, and Bruce~H. Edwards.
\newblock {\em Calculus: Early Transcendental Functions}.
\newblock Houghton Mifflin, fourth edition, 2007.

\bibitem{Marsden2004Vector-Calculus}
Jerrold~E. Marsden and Anthony~J. Tromba.
\newblock {\em Vector Calculus}.
\newblock W. H. Freeman, 5th edition, 2004.

\bibitem{Stewart2008Calculus:-Early}
James Stewart.
\newblock {\em Calculus: Early Transcendentals}, volume~6E of {\em Calculus}.
\newblock Brooks Cole, 2008.

\end{thebibliography}

\end{document}